\documentclass[12pt,psamsfonts]{amsart}
\usepackage{amsmath,amsthm,amsfonts,amssymb}
\usepackage{eucal,amscd}
\usepackage{pstricks}
\usepackage{bm}
\usepackage{epstopdf}
\usepackage{color}
\usepackage[foot]{amsaddr}
\usepackage{graphicx}
\usepackage[all,knot]{xy}
\xyoption{arc}

\newcommand{\ot}{\otimes}
\newcommand{\mC}{\mathcal{C}}
\newcommand{\Irr}{\textrm{Irr}}
\newcommand{\Z}{\mathbb{Z}}
\newcommand{\Q}{\mathbb{Q}}
\newcommand{\R}{\mathbb{R}}
\newcommand{\SU}{\textrm{SU}}
\newcommand{\M}{\mathcal{M}}

\begin{document}

\title{Quantum Computing: a Quantum Group Approach}

\author{Zhenghan Wang}
\address{Microsoft Station Q\\CNSI Bldg Rm 2237\\
    University of California\\
    Santa Barbara, CA 93106-6105\\
    U.S.A.\\Email: zhenghwa@microsoft.com}

\begin{abstract}

There is compelling theoretical evidence that quantum physics will change the face of information science.  Exciting progress has been made during the last two decades towards the building of a large scale quantum computer.  A quantum group approach stands out as a promising route to this holy grail, and provides hope that we may have quantum computers in our future.  

\end{abstract}

\maketitle

\section{Quantum group and quantum computing}

Working at Microsoft Station Q, I am interested in building a quantum computer by braiding non-abelian anyons \cite{Wang}.
The mathematical underpin of anyon theory is unitary modular tensor category theory---
an abstraction of the representation theory of quantum groups \cite{Kitaevhoneycomb} \cite{Wang}.
A quantum group commonly refers to a certain Hopf algebra such as a quasi-triangular one, though the terminology is used only loosely.
I would like to advocate the use of a quantum group to mean a modular category.  Modular categories are elaborate algebraic structures that arise in a variety of subjects: Hopf algebras, von Neumann algebras, topological quantum field theories (TQFTs), conformal field theories (CFTs), and topological phases of matter.   They are quantum analogues of finite abelian groups in a sense that I will explain below.  We could use unitary fusion categories to describe quantum symmetries of quantum phases of matter just as finite groups are used to describe symmetries of crystals.

"Qubits" are just the modern abstraction of a two-level quantum system and its composites.
Theoretically, there is compelling evidence that quantum physics can change the face of information science: the BB$84$ private key protocol\cite{BB84}, Shor's factoring algorithm\cite{Shor84}, and Hastings' additivity counterexamples\cite{Hastings09}.
But for the last two decades, the most interesting challenge has been to realize quantum computation in the real world.  A key "post-Shor" idea is to use topology to protect quantum information---a revolutionary idea articulated by A.~ Kitaev as follows:  \lq\lq if a physical system were to have quantum topological (necessarily nonlocal) degrees of freedom, which were insensitive to local probes, then information contained in them would be automatically protected against errors caused by local interactions with the environment.
This would be fault tolerance guaranteed by physics at the hardware level, with no further need for quantum error correction, i.e., topological protection."\cite{Freedman98}\cite{Kitaev97}

Materials with such topological degrees of freedom are called topological phases of matter\cite{Nayak08}.
The prospect of utilizing topological phases of matter to build a quantum
computer has been greatly enhanced by recent experimental progress.  This approach to quantum computing
at the triple juncture of mathematics, computer science, and physics is referred to as topological quantum computation (TQC).  The low energy effective theories of topological phases of matter are TQFTs, and their 
algebraic data are modular categories.  Therefore, unitary modular categories (UMCs) are
algebraic models of anyons in topological phases of matter.  TQC leads to new perspectives on TQFTs and modular categories.

\section{Modular categories as finite abelian groups}

A modular category $\mC$ is a non-degenerate ribbon fusion category \cite{ENO02}\cite{BK01}.  A ribbon fusion category is a fusion category with compatible braiding and pivotal structures.  A practical way to describe such a complicated structure $\mC$ is to use a representative set $\Irr(\mC)=\{X_i\}_{i\in I}$ of the isomorphism classes of simple objects of $\mC$, and the various structures on $\Irr(\mC)$.  The finite index set $I$ of $\Irr(\mC)$ is the isomorphism classes of simple objects of $\mC$, and an element of $I$ will be called a label.  The cardinality $|I|$ of the label set $I$ is called the {\it rank} of the modular category.  It is important to keep in mind the difference between labels and the simple object representatives.  To make an analogy with a finite group, we will think of each simple object $X_i$ as an analogue of a group element in the modular category $\mC$.  Then there will be a binary operation on $\Irr(\mC)$ which satisfies the analogues of associativity, commutativity, and inverse of a group.

The binary operation on $\Irr(\mC)$ is called the tensor product.  Just as tensor products of irreps of algebraic structures such as Hopf algebras, the tensor product $X_i\ot X_j$ is not necessarily an irrep, so it could decompose into a direct sum $\oplus_{k=1}^{|I|} N_{i,j}^k X_k$, i.e., $X_i\ot X_j \cong \oplus_{k=1}^{|I|} N_{i,j}^k X_k$ for some non-negative integers $N_{i,j}^k$ (called multiplicities).  This is some kind of quantum multiplication as the outcome is a superposition.  Actually the outcomes $\{X_k\}$ are endowed with an intrinsic probability distribution related to the quantum dimensions of the simple objects.

The associativities between $(X_i\ot X_j)\ot X_k$ and $X_i\ot (X_j\ot X_k)$ are weakened to be functorial isomorphisms that satisfy the pentagon axioms.

The dual $X_i^*$ of $X_i$ is quantum analogue of an inverse.  As such we would only require that the tensor unit $\bf{1}$ appears in $X_i\ot X_I^*$ as a sub-object.

Finally, the commutativity of $X_i\ot X_j$ and $X_j\ot X_i$ is relaxed to be a functorial isomorphism that satisfies the hexagon axioms.

There is a very convenient picture calculus for a modular category.  Basically, the axioms of a ribbon fusion category is to define topological invariants of colored braided trivalent framed graphs in the plane.  In particular, we have invariants of colored framed links.

Important numerical data for a modular category $\mC$ are the fusion multiplies $\{N_{i,j}^k\}, i,j,k\in I$, the $S$-matrix and $T$-matrix.  The fusion multiplies are organized into matrices $N_i, i\in I$ such that $(N_i)_{j,k}=N_{i,j}^k$.  The unnormalized $\tilde{S}=(S_{ij}$)-matrix consists of link invariants of the colored Hopf link: $S_{ij}$ is the link invariant of the Hopf link with the two components colored by $X_i$ and $X_j$. We always identify $I$ with $\{0,1,...,|I|-1\}$, and assume $X_0$ is the tensor unit $\bf{1}$.  The $T$-matrix $T$ is diagonal with the $i$-th diagonal entry equal to the twist $\theta_i$ of $X_i$.
The first row of the $\tilde{S}$ matrix consists of the invariants of the unknot colored by $X_i$, denoted as $d_i$, called the quantum dimension of the simple object $X_i$.  Let $D=\sqrt{\sum_{i\in I}d_i^2}$, $p_{\pm}=\sum_{i\in I}\theta_i^{\pm} d_i^2$.  Then an important identity for a modular category is:

$$p_+p_-=D^2.$$

The $S$-matrix $S$ is $\frac{1}{D}\tilde{S}$, and $S, T$ together leads to a projective representation of the modular group $SL(2,\mathbb{Z})$.

Any finite abelian group $G$ can be endowed with modular category structures in many different ways so that $G$ becomes the $\Irr(\mC)$'s.  One way to do it is as follows:  choose an $n\times n$ even integral matrix $B$ so that the cokernel of the map $B: \Z^n \rightarrow \Z^n$ is isomorphic to $G$, i.e., $\Z^n/\text{Im}(B)\cong G$.  Each element $g\in G$ can be identified as an equivalence class of vectors $v_g\in \Q^n$ such that $Bv_g\in \Z^n$.  Two such $v_g$'s are equivalent if their difference is in $\Z^n$.  Then there is a modular category $G_B$ with all $d_i=1$, fusion rules as $G$, and $\tilde{S}$-matrix ${S}_{ij}=e^{2 \pi \sqrt{-1} <v_{g_i}, v_{g_j}>_B}$, $\theta_{i}=e^{\pi \sqrt{-1} <v_{g_i}, v_{g_i}>_B}$, where $<v,w>_B=v^tBw$.

It is obvious that for a fixed order, there exist only finitely many groups.  The analogous statement for modular categories would be for a fixed rank, there are only finitely many equivalence classes of modular categories\cite{BNRW12}.  This rank finiteness property has been proved recently.  This fundamental theorem implies that modular categories are amenable to a classification, especially the low rank cases.  Modular categories up to rank $5$ have been classified\cite{BNRW12}\cite{RSW}.

\section{Modular categories as topological phases of matter}

Solid, liquid, and gas are all familiar states of matter.  But by a more refined classification, each state
consists of many different phases of matter.  For example, within the
crystalline solid state, there are many different crystals distinguished by
their different lattice structures.  All those phases are classical in the sense
they depend crucially on the temperature.  More mysterious and challenging to
understand are quantum phases of matter: phases of matter at zero temperature
(in reality very close to zero).  Quantum phases of matter whose low energy
physics can be modeled well by TQFTs are called topological phases of matter.
There are two classes of topological phases of matter in real materials: two
dimensional electron systems that exhibit the fractional quantum Hall effect
(FQHE), and topological insulators and superconductors.  The universal
properties of those materials can be modeled well by TQFTs, and their variations---symmetry enriched TQFTs.

Traditionally, phases of matter are classified by classical symmetries modeled by groups.
Quantum phases of matter cannot be fit into this Landau paradigm.  Intensive
effort in physics right now is to develop a post-Landau paradigm to classify
quantum phases of matter.  The well-studied classes of quantum phases of matter
are topological phases of matter.  To characterize topological phases of matter,
we focus on both the ground states and first excited states.  The
algebraic models of elementary excitations of the topological phases are unitary modular categories (UMCs).  If a classical symmetry modeled by a group $G$ is also
present, then the UMC is $G$-graded.  An interesting
example is the trivial case with symmetry including the topological insulators
discovered in the last few years---symmetry protect topological phases of matter.

Anyons are localized quantum fields.
The mathematical model of an anyon under UMCs is a simple object $X$ in a UMC, and in picture calculus, described as a small arrow. Therefore its world line in $\R^3$ is a ribbon.  The isomorphism class $x$ of an anyon $X$ is called the anyon type of $X$, or the anyonic charge.

A dictionary of terminologies between anyon theory and UMC theory is given in the following table.
\begin{table}[hb]
\begin{tabular}{|l|l|}\hline
\emph{UMC} & \emph{Anyonic system}\\\hline
simple object & anyon\\\hline
label & anyon type or anyonic charge\\\hline
tensor product & fusion\\\hline
fusion rules & fusion rules\\\hline
triangular space $V^c_{ab}$ or $V^{ab}_c$ & fusion/splitting space\\\hline
dual & antiparticle\\\hline
birth/death & creation/annihilation\\\hline
mapping class group representations & generalized anyon statistics\\\hline
nonzero vector in $V(Y)$ & ground state vector\\\hline
unitary $F$-matrices & recoupling rules\\\hline
twist $\theta_x = e^{2\pi i s_x}$ & topological spin\\\hline
morphism & physical process or operator\\\hline
colored braided framed trivalent graphs & anyon trajectories\\\hline
quantum invariants & topological amplitudes\\ \hline
\end{tabular}
\caption{}
\label{table-dictionary}
\end{table}

A UMC leads to an Atiyah-Segal type unitary TQFT\cite{Turaev94}.  An $(n+1)$-Atiyah-Segal-type TQFT is a symmetric monoidal functor from the $(n+1)$-bordism category $\textrm{Bord}(n+1)$ of $n$ and $n+1$-manifolds to the category of finitely dimensional vector spaces.  An important question is when such a TQFT has a Hamiltonian realization on triangulations such as the quantum double of finite groups and Turaev-Viro type $(2+1)$-TQFTs\cite{Kitaev97}\cite{LW05}.

\section{Modular categories as anyonic quantum computers}

Given $n$ anyons of type $x$ in a disk $D^2$, their ground state degeneracy
$\textrm{dim}(V(D,x,\cdots ,x))$ has the asymptotic growth rate $d_x$---the quantum dimension of $x$.
An anyon $X$ with $d_X=1$ is called an abelian anyon, and
an anyon with $d_X >1$ an non-abelian anyon.

Given a non-abelian anyon of type $x$, we can construct an anyonic quantum computation model as follows:
choose $V_{n,x}^a$ as the computational subspace---the Hilbert space of $n$ anyons of type $x$ and boundary
label $a$---the total anyonic charge.

To simulate a traditional quantum circuit
$U_L:({\mathbb{C}^2})^{\otimes n} \to ({\mathbb{C}^2})^{\otimes n} $ we need to find a braid $b$ making the square
\[\xymatrix{
({\mathbb{C}^2})^{\otimes n} \ar[r]^\iota \ar[d]_{U_L} &  V^a_{n,x} \ar[d]^{\rho(b)}\\
({\mathbb{C}^2})^{\otimes n} \ar[r]_\iota & V^a_{n,x}
}\]
commute, where $\rho(b)$ is the braid matrix and $\iota$ is an efficient embedding of the $n$-qubit space
$({\mathbb{C}^2})^{\otimes n}$ into the ground states $V^a_{n,x}$.
This is rarely achievable.
Therefore we seek $b$ making the square commute up to arbitrary precision.
To achieve universality for quantum
computation, we need to implement a universal gate set of unitary
matrices.  Then universality for anyonic quantum computing becomes the
question: can we find braids $b$ whose $\rho(b)$'s approximate a
universal gate set efficiently to arbitrary precision?  For
non-abelian anyons, since the Hilbert spaces always grow exponentially, universality is usually guaranteed if the
braid group representations afforded by the UMC have a dense
image in the special unitary groups $\SU(V^a_{n,x})$.

\subsubsection{Braiding universal anyons}

To simulate the standard quantum circuit model with braidings alone, we choose a non-abelian anyon type $x$ and use $m$ and $n$ anyons respectively to encode one and two qubits.
If the braid group representations $B_m \to \SU(V_{m,x}^a)$ and $B_n \to \SU(V_{n,x}^b)$ are independently dense for all anyon types $a,b$, then the anyonic quantum computing model is universal.  Anyon with such properties will be called braiding universal.  An example of a brading universal anyon is the Fibonacci anyon $\tau$.

\subsubsection{Resource assisted universal anyons}

It is conjectured that anyons with quantum dimensions $d_x^2 \notin \Z$ are essentially braiding universal.  But current experimentally accessible anyons are those $d_x^2\in \Z$.  An example of such anyons is the Ising anyon $\sigma$ with $d_{\sigma}=\sqrt{2}$.  More examples are the metaplectic anyons\cite{HNW12}. It is conjectured that such anyons always lead to representations of the braid groups with finite images, hence not braiding universal\cite{NR}.  Therefore, to obtain a universal braiding gate set, we need to supplement the braidings with other gates such as measurement, magic state distillation, ...  We will call an anyon $R$-assisted universal if the braidings supplemented with the resource $R$ give rise to a universal gate set.  For example, the Ising anyon is (magic-state-distillation)-assisted universal\cite{BK04}.

\section{Majorana physics and its realization}

At the writing, we are the closest to detecting and harnessing Majoranas.  Majoranas mean various things in the literature, and sometimes are very misleading.  We will distinguish between Majorana fermions from Majorana zero modes.  By Majorana fermions, we will mean abelian anyons with twist $\theta=-1$.  Therefore, Majorana fermions are never non-abelian.  An example of a Majorana fermion is the $\psi$ particle in the Ising theory.  By Majorana zero modes, we will mean the zero modes $\{\sum_\alpha c_\alpha |\psi_\alpha>\}$ obtained from the Majorana algebras $\mathcal{M}_n$ through the state-operator correspondence $\gamma=\sum_\alpha |\psi_\alpha> <\psi_\alpha|.$  In the Ising theory, the ground state manifolds of many Ising anyons $\sigma$'s consist of Majorana zero modes.  However, Majorana zero modes can arise without the presence of Ising anyons. 

The Majorana algebra $\mathcal{M}_n$ is generated by the operators $\{\gamma_i\}, i=1,\cdots, n-1$ such that $\gamma_i^{\dagger}=\gamma_i, \{\gamma_i, \gamma_j\}=2\delta_{ij}$. $\mathcal{M}_n$ supports a representation $\rho_{\mathcal{M}_n}$ of the braid group $B_n$ by sending the braid generator $\sigma_i$ to an automorphism $\rho_{\mathcal{M}_n}(\sigma_i)=\frac{1}{\sqrt{2}}(1-\gamma_i\gamma_{i+1})$ of the Majorana zero modes: $$\rho_{\mathcal{M}_n}(\sigma_i)\gamma = \frac{1}{\sqrt{2}}(1-\gamma_i\gamma_{i+1}) \cdot \gamma \cdot {\frac{1}{\sqrt{2}}(1-\gamma_i\gamma_{i+1})}^\dagger,$$
for any $\gamma \in \mathcal{M}_n$. 
 
It is easy to check:
$\begin{cases}
    \rho_{\M_n}(\sigma_i) {\rho_{\M_n}(\sigma_i)}^\dagger = \operatorname{Id},\\
    \rho_{\M_n}(\sigma_i)\gamma_i = \gamma_{i+1},\\
    \rho_{\M_n}(\sigma_i)\gamma_{i+1} = -\gamma_i,\\
    \rho_{\M_n}(\sigma_i)\gamma_j = \gamma_j, \text{ if $j\neq i$, $i+1$.}
\end{cases}$

This projective braid group representation  $\rho_{\mathcal{M}_n}(\sigma_i)=\frac{1}{\sqrt{2}}(1-\gamma_i\gamma_{i+1})$ is essentially the Jones representation of the braid groups at the $4$-th root of unity $q=e^{\frac{2\pi i}{4}}$ because the Temperley-Lieb Jones algebra at $q=e^{\frac{2\pi i}{4}}$ is the Majorana algebra $\mathcal{M}_n$\cite{Jones83}\cite{Wang}.

There are important experimental progress towards the confirmation of Majorana physics in real materials.  Interferometry experiments for the FQH liquids at filling fraction $\mu=\frac{5}{2}$ agree with the Ising physics\cite{Willett09}, and nanowires experiments are consistent with Majorana zero modes\cite{Leo12}. 

Recently, a generalization of Majorana physics to metaplectic physics is proposed\cite{HNW12}.  Potentially, metaplectic anyons can be used to perform universal quantum computation more easily than Majorana anyons.

\section{Open problems and future directions}

There are new problems in mathematics, physics, and computer science inspired by the quantum group approach to quantum computing.
The classification of modular categories is both an interesting question in mathematics and physics.  As beautiful mathematical structures, modular category theory could be regarded as quantum version of group theory.  In physics, a classification is very beneficial to the identification of topological phases of matter\cite{JWB12}.

As the theory of modular category reaches certain maturity, many new directions emerge.  UMCs capture topological properties of bosonic physical systems in $2$-spatial dimension.  Three directions that we would like to explore are for fermion systems, for boson systems with group symmetries, and for $3$-spatial dimension.

\subsubsection{SPT and SET}

Topological phases of matter can have conventional group symmetries $G$.  Such topological phases of matter are called symmetry enriched topological (SET) phases.  When the intrinsic topological order is trivial, SET is called symmetry protect topological (SPT) phase.  Important examples of SPTs are topological insulators and topological superconductors.  In $2$-spatial dimension, SETs are described by $G$-equivariant TQFTs/modular categories.  Classification of SETs then involves the interplay between group theory and modular category theory.   

\subsubsection{Fermion systems}

Fermion systems share similar Hilbert space as bosonic spin systems, but the meaning of locality is different\cite{GWW10}.  A local Hamiltonian for a fermion system is actually non-local.  Therefore, fermionic topological phases of matter are strictly richer than their bosonic counterparts.

\subsubsection{(3+1)-TQFTs}

Topological phases of matter behave differently in different spatial dimensions.  In $1D$, there are no intrinsic non-trivial topological order.  In $2D$, there is a richness of topological orders.  In $3D$, the situation is more intrigue.  $(3+1)$-TQFTs whose partitions that can detect smooth structures are very rare\cite{KW11}.  I would conjecture that all unitary Atiyah-Segal-type $(3+1)$-TQFTs give rise to partition functions that are simple homotopy invariants.  for $3D$, an algebraic framework analogous to modular category is lacking, at least we do not understand higher category for $3D$ to have a good understanding.

Quantum computers are still hypothetical models that promise revolutionary advance in information science.  The quantum algebra/topology approach provides substantial evidence that we may see quantum computers in our future, so the time of machines is coming.  

\section{Acknowledgments}The author is partially supported by NSF DMS 1108736.

\end{document}